\documentclass[11pt]{amsart}
\usepackage{amsmath,amsfonts,amsthm,amssymb,amscd}
\def\classification#1{\def\@class{#1}}
\classification{\null}

\DeclareFontFamily{OT1}{rsfs}{}
\DeclareFontShape{OT1}{rsfs}{n}{it}{<-> rsfs10}{}
\DeclareMathAlphabet{\mathscr}{OT1}{rsfs}{n}{it}

\newcommand{\F}{\mathbb{F}}

\newtheorem{lemma}{Lemma}

\newtheorem{theorem}{Theorem}
\theoremstyle{remark}
\newtheorem{remark}[theorem]{Remark}

\numberwithin{equation}{section}
\title{An improved sum-product inequality in fields of prime order}

\author{Misha Rudnev}
\address{Misha Rudnev, Department of Mathematics, University of Bristol,
  Bristol BS8 1TW, United Kingdom}
\email{m.rudnev@bristol.ac.uk}
\begin{document}
\begin{abstract} This note proves the sum-product inequality  $$\max(|A\cdot A|,\, |A+A|)\gg |A|^{1+\frac{1}{11}}\frac{1}{\log^\frac{4}{11}|A|},$$  \hspace{19mm}for sets of cardinality $|A|<\sqrt{p}$ in prime fields $\F_p$.
\end{abstract}
\maketitle
\section{Introduction}
Let $A$ be a subset of the multiplicative group $\F_p^*$ of the prime residue field $\F_p$ ($p$ is thought to be large, the constants throughout being independent of $p$). The sum set $A+A$ is defined as $\{a_1+a_2:\,a_1,a_2\in A\};$ in the same fashion one defines the difference, product, and ratio sets $A-A,\,A\cdot A,A:A$, respectively.

The sum-product estimate of Bourgain, Katz, and Tao (\cite{BKT}) and Konyagin (\cite{Ko}) establishes that as long as $|A|< p^{1-\Delta},$ one has
\begin{equation}
\max(|A+A|,\;|A\cdot A|)\;\geq c(\Delta) \; |A|^{1+\delta(\Delta)},
\label{spest}\end{equation}
for some $c(\Delta)$ and $\delta(\Delta)>0$.
In what follows, finite set cardinalities are denoted via $|\cdot|$; $c$ and $C$ denote some universal constants, changing from one place to another.

The nest question, along the lines of the Erd\"os-Szemer\'edi conjecture for the integers, is to establish quantitative estimates as to the relation $\delta(\Delta)$. Clearly, there is no no-trivial sum-product estimate in the case $\Delta=0$, i.e when the size of $A$ is commensurate with the size of the whole field. Therefore, in order to get a uniform quantitative estimate, one should impose a restriction on the size of $A$, which is usually $|A|<\sqrt{p}$, i.e $\Delta\geq \frac{1}{2}$, which is assumed further. (When $\Delta$ becomes smaller than $\frac{1}{2}$,  calculations involving exponential sums start yielding non-trivial results, rather quickly, as $\Delta$ decreases, overpowering the combinatorial estimates that this and the majority of the quoted papers deal with. See e.g. \cite{G1}. The case $|A|>\sqrt{p}$ is no longer discussed in the sequel.)

Garaev \cite{Ga} succeeded in obtaining the first quantitative sum-product estimate in $\F_p$: either the sum or the product set shall have cardinality $|A|^{1+\frac{1}{14}}$, up to a multiple of a power of $c\log |A|$. Katz and Shen (\cite{KS}) elaborated on a particular application of the Pl\"unnecke-Ruzsa inequality (cited here as Lemma \ref{pr}) in Garaev's proof and improved $\frac{1}{14}$ to $\frac{1}{13}$. Bourgain and Garaev (\cite{BG}) incorporated a covering argument (whose variant, taken from \cite{Sh}, is cited here as Lemma \ref{cover}) and improved $\frac{1}{13}$ to $\frac{1}{12}$. And Li (\cite{Li}) showed that a multiple of a power of $\log |A|$ can be done away with, so today's state-of-the art, for $|A|<\sqrt{p}$ is that
$$
\max(|A\cdot A|,\;|A+ A|)\;\geq \;c |A|^{1+\frac{1}{12}},
$$
for some absolute $c$. Moreover, one can replace either one or both the product set $A\cdot A$ with the ratio set $A: A$ and the sum set $A+ A$ with the difference set $A-  A$.
%

On the heuristic level, the proof of the above series of results effectively begins with assuming that the multiplicative energy of $A$ is large, which follows from the assumption that the product set (or the ratio set) is small. Then one chooses a ``popular'' dyadic interval, with approximately constant number of points of $A\times A$ on the lines through the origin in $\F_p\times\F_p$, corresponding to different ratios. The decision on the meaning of ``popular'' can be postponed to the end of the proof and depend on what the final estimates are. This fact that the notion of ``popularity'' could be adjusted to the final estimates, rather specified this in the ``natural'' way in the beginning of the proof was the observation due to Li (\cite{Li}), which enabled him to get rid of the factor $\log|A|$ which permeated the results preceding his.  Throughout the main body of the proof one uses the ``standard" arithmetic combinatorics technique: Pl\"unnecke-Ruzsa inequalities, covering lemmas and additive pivoting, which is roughly speaking estimating the size of $A+rA$ for some dilation $r$, in order to bring the sum or the difference set into play (as one wishes, since the choice of $A\pm A$ is allowed in both Pl\"unnecke-Ruzsa inequalities and the Covering lemma).

The technical issues involved in the estimates used to improve Garaev's original $\frac{1}{14}$ to $\frac{1}{12}$ have piled up, and this note offers yet another one, which enables to reduce ``the worst possible scenario''  estimating the size of $A+rA$ in the state-of-the-art literature, (that is this case ``to blame'' for the final answer $\frac{1}{12}$) to the ``normal case''  which would yield $\frac{1}{11}$, only if it were shown to be the only scenario.

This is done essentially via an application of the pigeonhole principle in Lemma \ref{focus} below which is used in the case study in the forthcoming proof.

\medskip
The main result is as follows.
\begin{theorem}\label{delta}
Let $A\subset \F^*_p$ with  $ |A| < \sqrt{p}$ and $p$ bigger than some absolute constant. Then, for some absolute $c$,  $$
\max(|A\cdot A|,\;|A+ A|)\;\geq \;c   |A|^{1+\frac{1}{11}}\frac{1}{\log^{\frac{4}{11}}|A|}.
$$
\begin{remark} One can replace either one or both the product set $A\cdot A$ with the ratio set $A: A$ -- in which case the logarithmic factor disappears -- and the sum set $A+ A$ with the difference set $A-  A$.\end{remark}

\end{theorem}
Notation-wise, the symbols $\ll,\,\gg,\approx$ suppress constants. E.g $|X|\gg|Y|$ means $|X|\geq c|Y|$, for some $c$.  The English language gets abused in accordance with these notations by saying ``at least'' or ``at most'' in the sense conveyed by the symbols $\gg,\ll$, respectively. To avoid confusion, ``at least'' and ``at most'' used in this way will be enclosed into quotation marks. E.g. saying that $|X|$ is ``at least'' $|Y|$ means $|X|\gg|Y|$.

\subsection{Acknowledgments}

The author would like to thank L. Li and A. Glibichuk for their comments on the initial draft of this paper, as well as T. Jones and O. Roche-Newton for helping him to identify errors in his earlier attempts to prove Theorem \ref{delta}.

\section{Lemmata}
$\;\;\;\;$ This section contains the two arithmetic combinatorics lemmas used.
The first one is a slight, but very useful generalisation of the Pl\"unnecke-Ruzsa inequality (see \cite{R}), which is due to Katz and Shen (\cite{KS}).
\begin{lemma}\label{pr}
Let $Y;\; X_1,\ldots X_k$ be additive sets. Then for any $0<\epsilon\leq 1$, there exists a subset $|Y'|\subseteq Y$, with $|Y'|\geq (1-\epsilon)|Y|,$ and some constant $C(\epsilon)$, such that
\begin{equation}\label{pr1}
|Y'+ X_1+\ldots+X_k| \;\leq  \;C(\epsilon)\,\frac{\prod_{i=1}^k|Y+X_i|}{|Y|^{k-1}}.
\end{equation}
\end{lemma}

The second one is a covering lemma, see e.g. \cite{BG}, \cite{Sh}, \cite{Li}.
\begin{lemma}\label{cover} Let $X_1$ and $X_2$ be additive sets. Then for any $\varepsilon\in (0,1)$ and some constant $C(\varepsilon)$, there exist ${\displaystyle \frac{C(\varepsilon)}{|X_2|}\min(|X_1+X_2|,|X_1-X_2|) }$ translates of $X_2$, whose union contains not less than $(1-\varepsilon)|X_1|$ elements of $X_1$.\end{lemma}

%
%


\section{Proof of Theorem \ref{delta}}
The proof goes essentially along the lines of the lately well-established approach, which in particular does not distinguish between the $A\pm A$, and applies as well to $A:A$, in which case it gets a little easier. If at the outset one deals with the product set, then the analysis begins by looking at the multiplicative energy whereupon it starts essentially dealing with ratios, rather than products. Throughout, one uses Lemmas \ref{pr} and \ref{cover} which admit any $\pm$ sign variation, and uses additive pivots $A+rA$, although it can be $A-rA$ instead. Hence, what follows is confined to the sum and the product sets only. Suppose both $|A+A|\leq K|A|$ and $|A\cdot A|\leq K|A|$, for some $K$.

First off, let us refine $A$, if necessary, so that the additive Pl\"unnecke-Ruzsa inequality applies to $A$ in full power as follows:

\begin{equation}\label{finepl}|A+A+A+A|\ll K^3|A|.\end{equation}
This will be used in the end of each of the Cases (i--iii) constituting the main body of the proof.

Indeed, if (\ref{finepl}) did not apply to $A$, one would, by Lemma \ref{pr}, choose a large subset $A'$ of $A$ (containing, say $90\%$ of its elements) such that with $A'$ replacing $A$, the estimate (\ref{finepl}) would be in place, and proceed with $A'$. As several more refinements will be made in the sequel, this first one is adopted without the change of notation: $A$ from now on stands for the possibly refined original set, which satisfies (\ref{finepl}) as well as the assumptions $|A+A|\leq K|A|$ and $|A\cdot A|\leq K|A|$.

The multiplicative energy of $E_*(A)$, i.e the number of ordered quadruples \\$(a_1,a_2,a_3,a_4)\in A\times A\times A\times A,$ satisfying the equation
$$
\frac{a_1}{a_3} = \frac{a_4}{a_2}
$$
is, by Cauchy-Schwarts, bounded from below as follows:
\begin{equation}\label{beg}
E_*(A)\geq \frac{|A|^4}{|A\cdot A|}\gg \frac{|A|^3}{K}.
\end{equation}
The multiplicative energy, by its definition, is
$$
E_*(A) = \sum_{\xi} n^2(\xi),
$$
where $\xi$ is the slope of a line through the origin in $\F_p^2$, and $n(\xi)$ is the number of points of $A\times A$ on this line.

Let us partition these lines, identified by their slopes, into dyadic groups as to how populated by points of $A\times A$ they are: a line in the $j$th group, $j=1,2,\ldots \leq 1+ \log_2 |A|$ will support $2^{j-1}\leq n(\xi)< 2^{j}$  points of $A\times A$. We shall consider a single ``popular'' dyadic group, to be specified in the end of the proof (although the choice will be straightforward). This dyadic group contains some number $L$ of lines, each supporting approximately $N$ points of $A\times A$, rendering the contribution $M$ into the multiplicative energy $E_*(A)$. One has
\begin{equation}\label{dy}
LN^2 \approx M.
\end{equation}
This implies, since $LN\ll |A|^2$, and $N\ll |A|$, that
\begin{equation}\label{lbn}
L,N\gg \frac{M}{|A|^2}.\end{equation}

Let $\Xi$ denote the set of slopes of the lines through the origin in the chosen dyadic group and $P$ denote the set of points of $A\times A$ supported on these lines. Since each line supports approximately the same number $N$ of points, one has
\begin{equation}
|P|\approx LN.\label{ptset}\end{equation}

Let \begin{equation} \begin{aligned}A_x&=&\{y : (x,y)\in P\},\\A_y&=&\{x : (x,y)\in P\}. \label{ai}\end{aligned}\end{equation}
Namely $A_x$ is the set of {\em ordinates} of points of $P$ all having the same abscissa $x$ and $A_y$ the set of {\em abscissae} of points of $P$ all having the same ordinate $y$. For $y\in A_x$ let $A_{y/x}\equiv A_\xi$ (with $\xi =\frac{y}{x}$) be the set of {\em abscissae} of the points of $P$ supported on the line through the origin with the slope $\xi =\frac{y}{x}$ passing through $(x,y)$. So, for all $y\in A_x$, $|A_{y/x}|\approx N.$

The next lemma is an important building block for the proof.
\begin{lemma}\label{focus}
There exist a popular abscissa $x=\tilde x$ and a popular ordinate $y=\tilde y$, i.e. such that $|A_{\tilde x}|, |A_{\tilde y}|\gg \frac{LN}{|A|}$, as well as a subset $\tilde A_{\tilde x}\subseteq A_{\tilde x}$, with
\begin{equation}
|\tilde A_{\tilde x}|\gg\frac{LM}{|A|^3},
\label{onebug}\end{equation}
such that for every $y\in \tilde A_{\tilde x}$, one has \begin{equation}\label{cl}
|\tilde A_{y}\equiv A_{y/\tilde x} \cap A_{\tilde y}|\gg \frac{LMN}{|A|^4}.\end{equation}
\end{lemma}

\subsubsection*{Proof of Lemma \ref{focus}}
To be precise, one must run some popularity arguments about the sets $P$ and $\Xi$. In essence however, since $|P|\approx LN$ and the maximum population of a line with slope in $\Xi$ is ``at most'' $ N$, this means that lower bounds pertaining to sums over the sets $\Xi, P$ also apply to the whole hierarchy of their ``popular subsets'', for merely the price of the constants hidden in the $\gg$ symbols getting worse.

\medskip
By the pigeonhole principle, a positive proportion of the points of $P$ have ``popular ordinates'', namely such that a horizontal line passing through such an ordinate contains ``at least'' $\frac{LN}{|A|}$ points of $P$. Let $A'$ denote the set of these popular ordinates.

Let us call $P'$ the set of the points of $P$ having popular ordinates. One still has $|P'|\gg LN$ and $P'$ is still supported on no more than $L$ lines passing through the origin. Therefore, $P'$ still carries a positive proportion of the multiplicative energy $M\approx LN^2$ of $P$, and an average line (identified by its slope in $\xi\in\Xi$) still has $\sim N$ points of $P'$.

Now the abscissae of the points in the set $P'$ can be pruned in the same way, down to the ``popular'' ones, namely such that there are ``at least'' $\frac{LN}{|A|}$ points of $P'$ supported on a vertical line with a popular abscissa, whereupon $P'$ gets refined to $P''$. One can now disregard the lines  that may have become ``poor'', namely those supporting fewer than $cN$ points of $P''$. The remaining lines have a set of slopes $\Xi'$, whose cardinality is still ``at least'' $L$, and each of these lines still supports $\approx N$ points of $P''$.

The claims of the Lemma are basically the average case for the sum
\begin{equation}
\Sigma = \sum_{x\in A'', y\in A'} \;\sum_{z\in A_x} \;\; |A_{z/x} \cap {A_{y}}|,
\label{sumsg}\end{equation}
Let us write $z\sim x$ if $(x,z)$ is a point of $P$.
Then
$$
\Sigma \gg \sum_{x\in A'',\,y\in A',\,z:\,z\sim x} |A_{z/x} \cap {A_{y}}|.
$$

Rearranging the summations by fixing the ratio $\xi = \frac{z}{x}$ and confining the sum only to the ``popular'' ratios in $\Xi'$ (these ratios each have ``at least'' $N$ realisations in terms of ``popular'' ordinates and abscissae) one gets:
\begin{equation}
\Sigma \gg  N \sum_{\xi\in \Xi'} \sum_{y\in A'} |A''_{\xi} \cap A_{y}|,
\label{sumsg1}\end{equation}
The sum in the right-hand-side of (\ref{sumsg}) takes place over ``at least'' $L$ lines whose slopes are popular, with respect to the refined point set $P''$, i.e. each line supporting ``at least'' $N$ points of the latter point set.
Fixing a particular line $l$, let $A_l$ be the set of the (popular) abscissae of the points of $P''$ on this line. One has $|A_l|\gg N$, as well as $A_l\subseteq A'',$ the set of popular abscissae. This means, by construction of the sets $A', A''$, that each element of $A_l$ belongs to ``at least'' $\frac{LN}{|A|}$ different sets $A_{y}$, with $y\in A'$.

Hence,

\begin{equation}\Sigma\gg N\cdot L\cdot N\cdot \frac{LN}{|A|}\gg \frac{LMN}{|A|}\label{sbd}.\end{equation}

Now  the claims (\ref{cl}) and (\ref{onebug}) follow  from (\ref{sumsg}--\ref{sbd}) by the pigeonhole principle. Indeed, there is a pair of $(\tilde x,\tilde y)\in A''\times A'$, such that
\begin{equation}\sum_{z\in A_{\tilde x}} \;\; |A_{z/\tilde x} \cap {A_{\tilde y}}|\gg \frac{LMN}{|A|^3}.\label{sbd1}\end{equation} One can restrict the summation in (\ref{sbd1}) to those values of $z$ only, for which
$$|A_{z/\tilde x} \cap {A_{\tilde y}}|\gg \frac{LMN}{|A|^4}.$$
These $z$ are to form the set $\tilde A_{\tilde x},$ whose existence is stated in the Lemma.
Since the maximum size of a single $A_{z/\tilde x}$ is ``at most'' $N$, then the number of such $z$ is ``at least'' $\frac{LM}{|A|^3}.\;\;\;\Box$

\bigskip
$A$ can be now scaled, so without loss of generality assume $\tilde x=1$, and to save on the notations let $A_{\tilde x}\equiv B$ (and scale no more). Let us also us further use the notations $C$ for the set $A_{\tilde y}$
and $\tilde B$ for the refinement of $B$, so that (\ref{cl}) holds for {\em every} $z\in \tilde B$. Heuristically, $C$ catches a large proportion of the abscissae of points of $P$ supported on every line through the origin, whose slope is in $\tilde B$.

With the parameter $r\in\F^*_p$, consider now (as a case study) ``additive pivots'' $C + r C$, $\tilde B+r\tilde B$, as well as finally $C +r \tilde C_p$, where $\tilde C_p$ will be subset of $C$ containing ``many'' abscissae of points of $P$ lying on the particular line to be specified, with the slope $p\in\tilde B$.

Introduce the following notation: for a pair of sets $S_1,S_2\subset \F_p$, let
\begin{equation}\label{er}
R(S_1,S_2) = \left\{\frac{p-q}{s-t}:\;p,q,\in S_1,\,s,t\in S_2, \;p\neq q,s\neq t\right\}.\end{equation}

One now faces with three cases that deal with the sets $R(C,C),\,R(\tilde B,\tilde B), $ and $R(C,\tilde C_p)$, where $\tilde C_p$ is the above-mentioned subset of $C$, to be chosen. Case (i) and (ii) will constitute the {\em normal case}, where the estimates somewhat repeat what has been done in the state-of-the-art literature. These estimates, from \cite{Ga} on, it would give a better estimate securing the exponent $\frac{1}{11}$, rather than $\frac{1}{12}$ with the current state of the art. But there would always be the ``worst possible scenario'' which forced the denominator in the final answer to become by one bigger than rendered by the normal case. This ``worst case scenario'', however, gets replaced below by Case (iii) which gives a similar estimate, as in the normal case.

Recall, by Lemma (\ref{focus}),
\begin{equation}\label{rest}  |\tilde B|\gg \frac{LM}{|A|^3},\qquad |C|\gg \frac{LN}{|A|}.\end{equation}

\medskip

{\sf Case (i):} Suppose $R(\tilde B,\tilde B)\neq R(C,C)$. Then there are two possibilities:

\medskip
$\;\;\;\;${\sf Case (i.1):} There is $r=\frac{p-q}{s-t}\in R(\tilde B,\tilde B),$ such that $r$ is not in $R(C,C).$ Fix this $r$ along with a quartet $p,q,s,t$ representing it throughout the rest of Case (i.1).

Then the equation
\begin{equation}\label{quad1}
a_1+ ra_2 = a_3 + r a_4,
\end{equation}
where $a_1,\ldots, a_4\in C,$ has only trivial solutions, i.e those with $a_2=a_4$.
Thus for any $C'\subseteq C$, which contains a positive proportion of the elements of $C$, one has
\begin{equation}\label{bchain1}
|C'+rC'| \gg |C|^2.
\end{equation}

Therefore,
\begin{equation}\label{vot}|C|^2 \ll \left|C' + \frac{p-q}{s-t}C'\right|\ll |p C' - qC' + sC' -tC'|.\end{equation}
But $p, q,s,t$ are elements of $\tilde B$, and hence of the set of slopes $\Xi$. Therefore one can bound the sizes of  the corresponding dilates of $C'$ from above, using Lemma \ref{cover}.

Lemma \ref{cover} specifies the refinement $C'$ of $C$ which is chosen in such a way that the sets involved get covered by translates of $A$ not merely in their large proportion as guaranteed by the Covering lemma, but by $100\%$.

Recall the notation $A_\xi\subseteq A$ for the set of the abscissae of the line with the slope $\xi.$ For any $\xi\in \Xi$, $99\%$ of the set $\xi C$ can be covered by ``at most''
\begin{equation}
\frac{|\xi C + \xi A_\xi|}{|\xi A_\xi|}  \ll \frac{|A+A|}{N} \ll \frac{K|A|}{N}
\label{cov}\end{equation}
translates of $\xi A_\xi\subseteq A$ (which have cardinality approximately $N$), and therefore translates of $A$ itself. Now one can choose $C'\subset C$, appearing in (\ref{bchain1}) as a subset containing at least every second element of $C$, and such that the sets $ p C',$ $ - q C',$ $sC',$ $- tC',$ each get fully covered (independently) by ``at most'' $\frac{K|A|}{N}$ translates of $A$ each. Which means, the whole set $(s-t) C'+ (p-q) C'$ has been fully covered by ``at most'' $\left(\frac{K|A|}{N}\right)^4$ translates of  $A+A+A+A$.

Together with (\ref{bchain1}) this implies
\begin{equation}
N^4 |C|^2\ll K^{4}|A|^4|A+A+A+A|\ll K^7|A|^5,
\label{fin1}\end{equation}
after
using (\ref{finepl}).

\medskip
$\;\;\;\;${\sf Case (i.2):} There is $r=\frac{p-q}{s-t}\in R(C,C),$ such that $r$ is not in $R(\tilde B,\tilde B).$ Fix this $r$ along with a quartet $p,q,s,t$ representing it. Recall that $C=A_{\tilde y}$, the set of abscissae with a popular, with regard to $P$, ordinate $\tilde y.$ Thus $\frac{\tilde y}{p},$ $\frac{\tilde y}{q}$ ,$\frac{\tilde y}{s},$ $\frac{\tilde y}{t}$ are the slopes in $\Xi$. And so are, since the point set $P$ is symmetric with respect to the bisector $y=x$, the ratios $\xi_p =  \frac{p}{\tilde y},$ $\xi_q =  \frac{q}{\tilde y},$ $\xi_s =  \frac{s}{\tilde y},$ and $\xi_t =  \frac{t}{\tilde y}.$

Thus
\begin{equation}
r=\frac{\xi_p-\xi_q}{\xi_s-\xi_t}.
\label{inv}\end{equation}
Repeat now the argument of Case (i.1) from (\ref{bchain1}) through (\ref{fin1}) replacing there $C$ with $\tilde B$ and $(p,q,s,t)$  with $(\xi_p,\xi_q,\xi_s,\xi_t)$. This yields the analogue of the estimate (\ref{fin1}):

\begin{equation}
N^4 |\tilde B|^2\ll K^{4}|A|^4|A+A+A+A|\ll K^7|A|^5.\label{fin11}\end{equation}

The estimate (\ref{fin11}) is slightly worse than (\ref{fin1}), since the estimate on the size of $\tilde B$ in Lemma \ref{focus} is more restrictive than on the size of $C$, see (\ref{rest}). So the estimate (\ref{fin1}) will not be returned to.

By Lemma \ref{focus}, one has $|\tilde B|\gg \frac{LM}{|A|^3}$, and as $M\approx LN^2$, (\ref{fin11}) implies that
\begin{equation}
M^4 \ll K^{7}|A|^{11}.
\label{fin2}\end{equation}

\medskip
From now on let us assume that $R(C,C)=R(\tilde B, \tilde B)\equiv R$.

\medskip

{\sf Case (ii):} Suppose
$|R|\gg \min(|C|^2,p)$. The original  assumption $|A|<\sqrt{p}$ ensures that the minimum is achieved when $|R|\gg |C|^2$. For $r\in R$, let $E_r(C)$ denote the number of ordered quadruples $(a_1,a_2,a_3,a_4)\in C\times C\times C\times C$ satisfying the equation (\ref{quad1}).

The equation (\ref{quad1}) has trivial solutions -- when $a_2=a_4$ -- and non-trivial ones -- when $(a_1,a_2,a_3,a_4)$ determines $r$. Summing over all $r$, one has, under the Case (ii) assumption:

\begin{equation}\label{ph}\sum_{r\in R}  E_r(C) = |R| |C|^2 + |C|^4\ll  |R| |C|^2 .\end{equation}
Hence, there exists $r\in R$, such that $E_r(C)\ll |C|^2$. Fix this $r$ throughout the rest of Case (ii), together with some quartet $p,q,s,t\in C$ representing it. Observe that the same bound is certainly satisfied by $E_r(C')$, namely when the equation (\ref{quad1}) is restricted to a subset $C'\subseteq C$. If $C'$ contains a positive proportion of the elements of $B$, then by Cauchy-Schwartz,
\begin{equation}\label{bchain2}
|C'+rC'| \gg |C|^2.
\end{equation}
Now use the representation (\ref{inv}) for $r$ and
 repeat verbatim -- with $(\xi_p,\xi_q,\xi_s,\xi_t)$ replacing $(p,q,s,t)$ there -- the argument from (\ref{bchain1}) on within Case (i.1).

\medskip
From now on assume that $R(C,C)=R(\tilde B, \tilde B)\equiv R$ and $|R|\ll |C|^2.$ This takes one into the final Case (iii) -- the analog of worst possible scenario in the former literature which turns out to be a regular one here. This, together with Lemma \ref{focus}, is central for the proof.

\medskip
{\sf Case (iii):}
 There exists some $r=\frac{p-q}{s-t}\in R$, such that $r+1\not\in R$ (otherwise $R=\F_p^*$, and the condition $|A|<\sqrt{p},$ returns one into Case (ii)). Fix this $r$, together with a quartet $p,q,s,t\in \tilde B$ representing it.

Consider now the quantity $p\in \tilde B$ and the set $\tilde C_p$, a subset of $C$, furnished by Lemma \ref{focus}, with ``at least'' $\frac{LMN}{|A|^4}$ elements (set in (\ref{cl}) $y=p$ and recall that $\tilde x = 1$ after scaling). The members of $\tilde C_p$ are the {\em abscissae} of points of $P$, supported on the line through the origin with the slope $p$. They also lie in the set $C$. By the assumption  of Case (iii) one has $R(C,C)=R(\tilde B, \tilde B)=R$, and it follows that since $r+1$ is not in $R$, it is not in $R(C,\tilde C_p)\subseteq R$.

Thus, using Lemma \ref{focus}, which guarantees that $|C|\gg \frac{LN}{|A|}$, as well as (\ref{cl}), for any positive proportion subsets $C',\tilde C_p' $ of  $C,\tilde C_p$, respectively, one has
\begin{equation}\label{bchain3}
|C' + (r+1) \tilde C'_p| \gg \frac{LN}{|A|}\frac{LMN}{|A|^4}.\end{equation}

Let us use Lemma \ref{pr} in the following form: for any positive proportion subset $C'$ of $C$ ($C'$ is to be specified later, in the ensuing covering argument), there is $C''\subseteq C'$, containing, say $99\%$ of the members of $C'$, and such that, by (\ref{bchain3}), one has:
\begin{equation}\begin{aligned} \frac{LN}{|A|}\frac{LMN}{|A|^4} \ll |C'' + (r+1) \tilde C_p'| &\ll \left|C'' + \tilde C'_p + \frac{p-q}{s-t} \tilde C_p'\right|   \\ &\ll
\frac{|C'+ \tilde C'_p|}{|C'|} \left|C' + \frac{p-q}{s-t} \tilde C_p'\right|
 \\ &\ll
  \frac{K|A|}{|C|}\left|(s-t)C' + (p-q)\tilde C_p'\right|
  \\ &\ll
  \frac{K|A|}{|C|}\left|sC'-t C' -q\tilde C_p' + A\right|.\end{aligned}\label{nda}\end{equation}
Indeed, $p \tilde C_p \subseteq A$, since $p\in \Xi$ is a slope and $\tilde C_p$ is a subset of the abscissae of the points of $P$ supported on the corresponding line through the origin. In comparison to the above Cases (i,ii), the covering argument is now to be applied only three rather than four times.

The covering argument now proceeds in exactly the same way as it was in Cases (i,ii), since $s,t,q$ are members of $B\subseteq \Xi$. The subsets $C'$ of $C$, as well as $\tilde C_p'$ of $\tilde C_p$ are chosen so that
their dilates $sC'$,  $-tC'$,  $-q\tilde C_p'$ are each fully covered (independently) by ``at most'' $\frac{K|A|}{N}$ translates of $A$. (E.g. $-q\tilde C_p'$ gets covered by no more than
$\frac{|-q\tilde C_p' - q A_q|}{|A_q|}\ll \frac{K|A|}{N}$ translates of $q A_{q}\subseteq A$, and therefore $A$ itself.)

Then (\ref{nda}) yields
$$
\frac{LMN}{|A|^4}\left(\frac{LN}{|A|}\right)^2 N^3\ll K^4|A|^4|A+A+A+A|\ll K^7|A|^5,
$$by (\ref{finepl}). Since $M\approx LN^2$, this means:
$$
M^4\ll K^7|A|^{11},
$$
restating the estimate (\ref{fin2}), which therefore
shall be taken as the final one in the above Case study.

\medskip
To conclude the proof, all one has to do is choose $M$. In view  of the latter estimate, which only contains the multiplicative energy $M$, one can only do the straightforward pigeonholing
\begin{equation}\label{choiceofm}
M\gg \frac{E_*(A)}{\log|A|}\gg \frac{|A|^3}{K\log|A|}, \end{equation} since by the pigeonhole principle, there is a dyadic group of lines through the origin furnishing ``at least'' the above amount of multiplicative energy.

Therefore
\begin{equation}\label{withlog}
K\gg \left(\frac{|A|}{\log^4|A|}\right)^{\frac{1}{11}}.\end{equation}


This completes the proof of Theorem \ref{delta}. $\Box$

\medskip
The final remark is that having in the outset the ratio set $A: A$, rather than the product set $A\cdot A$ makes the multiplicative energy argument superfluous and matters more transparent. $\Xi$ would be denote the set of ``popular lines'' with the slopes in $\Xi$, each supporting ``at least'' $\frac{|A|}{K}$ points of $A\times A$. This rids one of the necessity of using dyadic pigeonholing with the set of slopes and as the result, there is no logarithm in the final estimate. Altogether the lines with the slopes in $\Xi$ support a positive proportion of the set $A\times A$. The sets  $B, \tilde B,$ and $C$ then become positive proportion subsets of $A$, and the relations affecting the left-hand side of the estimate (\ref{fin2}) would assume the simpler form $LN\gg |A|^2,\,N\gg\frac{|A|}{K}$.


\bigskip
\bigskip
\bigskip

\begin{thebibliography}{4}

\bibitem{BG} J. Bourgain, M.Z. Garaev. {\em On a variant of sum-product estimates and explicit exponential sum bounds in prime fields.}
Math. Proc. Cambridge Philos. Soc. 146 (2009), no. 1, 1--21.

\bibitem{BKT} J. Bourgain, N. Katz and T. Tao. {\em A sum-product estimate in finite fields and their applications.} Geom. Func. Anal. 14 (2004), 27–-57.

\bibitem{Ga} M.Z. Garaev. {\em An explicit sum-product estimate in} $\F_p$. Intern. Math. Res. Notices (2007), no 11, 1--11.

\bibitem{G1} M.Z. Garaev. {\em The sum-product estimates for large subsets of prime fields.} Proc. Amer.
Math. Soc. 137 (2008), 2735 -- 2739.



\bibitem{Ko} S.V. Konyagin. {\em A sum-product estimate in fields of prime order.} Preprint {\sf arXiv:math/0304217} (2003), 9pp.

\bibitem{KS} N.H. Katz, C.-Y. Shen.  {\em A slight improvement to Garaev's sum product estimate.} Proc. Amer. Math. Soc. 136 (2008), 2499--2504.

\bibitem{Li} L. Li. {\em Slightly improved sum-product estimates in fields of prime order.} Preprint {\sf arXiv:math/0907.2051} (2009), 9pp.

\bibitem{R}  I. Z. Ruzsa.  {\em An application of graph theory to additive number theory.} Scientia, Ser. A 3 (1989), 97--109.


\bibitem{Sh} Ch.-Y. Shen. {\em An extension of Bourgain and Garaev's sum-product estimates.} Acta Arith. 135 (2008), 351–-356.

\end{thebibliography}
\end{document}